\input{aipcheck}

\documentclass[final]{aipproc}
\layoutstyle{8x11single}

\usepackage{autograph,epic}
\usepackage{latexsym,bezier,amsbsy,psfrag,color,enumerate,amsfonts,amsmath,amscd,amssymb,amsthm}
\usepackage{epsfig}
\setloopdiam{17}\setprofcurve{10}

\newtheorem{definition}{Definition}
\newtheorem{theorem}[definition]{Theorem}

\newtheorem{example}[definition]{Example}

\begin{document}

\title{Computing the Wiener index in Sierpi\'{n}ski carpet graphs}


\author{Daniele D'Angeli}{
address={Institut f\"{u}r Mathematische Strukturtheorie (Math C), TUGraz, Steyrergasse 30, 8010 Graz, Austria. \footnote{Email address: dangeli@math.tugraz.at \quad (The author was supported by Austrian Science Fund project FWF P24028-N18.)}}  }

\author{Alfredo Donno}{
address={Universit\`{a} degli Studi Niccol\`{o} Cusano - Via Don Carlo Gnocchi, 3 00166 Roma, Italia. \footnote{Email address: alfredo.donno@unicusano.it \quad (Corresponding author)}}}

\author{Alessio Monti}{
address={Universit\`{a} degli Studi Niccol\`{o} Cusano - Via Don Carlo Gnocchi, 3 00166 Roma, Italia. \footnote{Email address: alessio.monti@unicusano.it}}}

\begin{abstract}
We describe an algorithm to compute the Wiener index of a sequence of finite graphs approximating the Sierpi\'{n}ski carpet.
\end{abstract}

\maketitle


\section{Introduction}
The famous Sierpi\'{n}ski carpet is a self-similar, infinitely ramified fractal introduced by
W. Sierpi\'{n}ski in 1916 \cite{courbe}. Many physical models and critical phenomena on this and other related fractals, or on their discrete approximations, have been widely studied in the literature \cite{bonnier, dimeri, ising, gefen2, gefen3, percolation}. In this short note, we study the Wiener index of a sequence $\{\Gamma_n\}_{n\geq 1}$ of graphs studied in \cite{rodi, compactification}, which forms a discrete approximation of the Sierpi\'{n}ski carpet. The \textit{Wiener index} $W(\Gamma)$ of a graph $\Gamma$ concerns the sum of distances between vertices in a finite graph. It was introduced by the chemist H. Wiener \cite{Wiener} in order to find correlations between physicochemical properties of organic compounds, and the topological structure of their molecular
graphs. This index is one of the most studied topological invariants in
mathematical chemistry. It has been investigated for many classes of graphs (see, for instance, \cite{sandi1, sandi2, donnoiacono}). We construct an explicit algorithm to compute the value $W(\Gamma_n)$, by using a suitable embedding of the graphs $\Gamma_n$ into $\mathbb{Z}^2$.

\section{Wiener index and Carpet graphs}

Let $\Gamma=(V,E)$ be a connected finite graph. The \textit{Wiener index} $W(\Gamma)$ of $\Gamma$ is the number
$$
W(\Gamma)=\frac{1}{2}\cdot\sum_{v,w\in V} d(v,w),
$$
where $d$ is the geodesic distance on $\Gamma$. We investigate it for the sequence of \textit{Carpet graphs} $\{\Gamma_n\}_{n\geq 1}$ defined below.

\begin{figure}[h]
\psfrag{3}{$3$}\psfrag{4}{$4$}\psfrag{5}{$5$}\psfrag{6}{$6$}\psfrag{7}{$7$} \psfrag{0}{$0$}\psfrag{1}{$1$}\psfrag{2}{$2$}
 \psfrag{a}{$a$}\psfrag{b}{$b$}\psfrag{c}{$c$}\psfrag{d}{$d$}
 \psfrag{dd}{$c$}\psfrag{d4}{$c4$}\psfrag{d45}{$c45$}
 \psfrag{A3}{$A_3$}    \psfrag{B3}{$B_3$} \psfrag{C3}{$C_3$} \psfrag{D3}{$D_3$} \psfrag{Gammabarra}{$\overline{\Gamma}$}
 \psfrag{C4}{$C_4$}    \psfrag{Gamma1}{$\Gamma_1$} \psfrag{Gamma2}{$\Gamma_2$} \psfrag{Gamma3}{$\Gamma_3$} \includegraphics[width=0.7\textwidth]{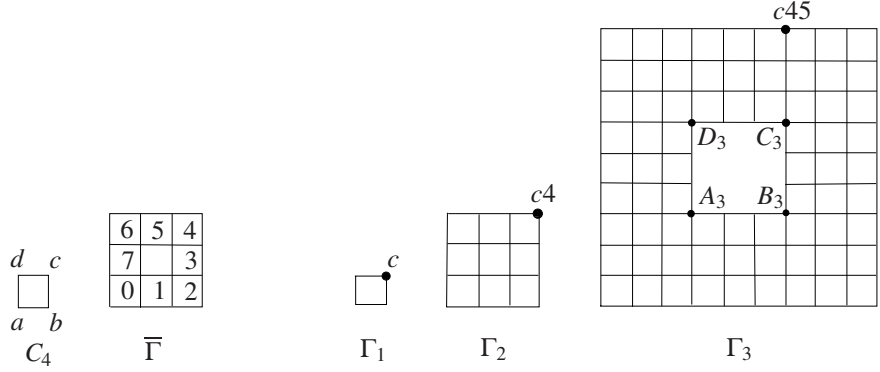}
\caption{Construction of the graph $\Gamma_n$ and examples.}\label{fig1}
\end{figure}

Fix two finite alphabets $X=\{0,1,\ldots, 7\}$ and $Y=\{a,b,c,d\}$, and let $X^\infty=\{x_1x_2\ldots : x_i\in X\}$ be the set of all infinite words over $X$. Let $C_4$ be the cyclic graph of length $4$ whose vertices will be denoted by $a, b, c, d$ (Fig. \ref{fig1}). \\

\textit{\textbf{Recursive construction of the graphs $\Gamma_n$}}.
\begin{itemize}
  \item \textbf{Step $1$.} The graph $\Gamma_1$ is the cyclic graph $C_4$. 
  \item \textbf{Step $n-1 \to n$.} Take $8$ copies of $\Gamma_{n-1}$ and glue them
  together on the model graph $\overline{\Gamma}$, in such a way that these copies occupy the positions indexed by $0,1,\ldots, 7$ in $\overline{\Gamma}$ (Fig. \ref{fig1}). Note that each copy shares at most one (extremal) side with
  any other copy. A vertex of $\Gamma_n$ is associated with the word $yx_1\ldots x_{n-1}\in Y\times X^{n-1}$, if it belongs to the copy of $\Gamma_{n-1}$ indexed by $x_{n-1}$ contained in $\Gamma_n$, and if it was associated with the word $yx_1\ldots x_{n-2}$ at the previous step. Notice that any word in $Y\times X^{n-1}$ corresponds to a unique vertex of $\Gamma_n$, but the viceversa is not true.
\end{itemize}
 In Fig. \ref{fig1}, for instance, the vertex associated with the word $c45$ in $\Gamma_3$ also corresponds to the word $d64$.

\section{Distances and numerical results}

The graph $\Gamma_n$ can be embedded into $\mathbb{Z}^2$, in such a way that the left corner vertex on the bottom of $\Gamma_n$ coincides with the vertex $(0,0)$ of $\mathbb{Z}^2$, and each horizontal edge of $\Gamma_n$ coincides with an edge of $\mathbb{Z}^2$ connecting two vertices of type ${\bf v}$ and ${\bf v}\pm (1,0)$, whereas each vertical edge of $\Gamma_n$ coincides with an edge of $\mathbb{Z}^2$ connecting two vertices of type ${\bf v}$ and ${\bf v}\pm (0,1)$. This embedding does not preserve, in general, the distances. In fact, the bigger and bigger holes, denoted by $H_n$, that we produce in our recursive construction may contain edges that belong to the shortest path connecting two vertices of $\Gamma_n$ after its embedding into $\mathbb{Z}^2$. In some sense, these holes correspond to \textit{obstructions} (in the terminology of \cite{compactification}) that should be taken into account in order to determine $W(\Gamma_n)$. In what follows, we analyse all possible cases in which the distances between vertices differ from the case $\mathbb{Z}^2$: we will say that the path joining such vertices meets an obstruction. The representation of the vertices of $\Gamma_n$ in terms of words in $Y\times X^{n-1}$, the description of the obstructions and the symmetry of the graphs $\Gamma_n$ allow us to describe an algorithm to compute $W(\Gamma_n)$.

In order to explicitly describe the embedding of $\Gamma_n$ into $\mathbb{Z}^2$, we define the following vectors of $\mathbb{Z}^2$ associated with the letters of the alphabet $Y$:
$$
{\bf v}_a = (0,0) \quad  {\bf v}_b = (1,0) \quad {\bf v}_c = (1,1) \quad{\bf v}_d = (0,1)
$$
and the following vectors of $\mathbb{Z}^2$ associated with the letters of the alphabet $X$:
$$
{\bf v}_0= (0,0) \quad  {\bf v}_1 = (1,0) \quad {\bf v}_2 = (2,0) \quad{\bf v}_3 = (2,1) \quad
{\bf v}_4 = (2,2) \quad  {\bf v}_5 = (1,2) \quad {\bf v}_6 = (0,2) \quad{\bf v}_7 = (0,1).
$$
Now let $w_n = yx_1x_2\ldots x_{n-1} \in Y\times X^{n-1}$; we associate with $w_n$ a vector ${\bf w}_n$ of $\mathbb{Z}^2$ defined as ${\bf w}_n = {\bf v}_y+ \sum_{i=1}^{n-1}3^{i-1}{\bf v}_{x_i}$. Given two finite words $w_n^1$ and $w_n^2$ and the corresponding vectors ${\bf w}_n^1(X({\bf w}_n^1),Y({\bf w}_n^1))$ and ${\bf w}_n^2(X({\bf w}_n^2),Y({\bf w}_n^2))$, we put $\|{\bf w}_n^1-{\bf w}_n^2\|_1 = |X({\bf w}_n^1)-X({\bf w}_n^2)|+|Y({\bf w}_n^1)-Y({\bf w}_n^2)|$, which is the geodesic distance in $\mathbb{Z}^2$.
Given $w_n^1= y^1x_1^1\ldots x_{n-1}^1, w_n^2 =y^2x_1^2\ldots x_{n-1}^2\in Y\times X^{n-1}$, let $h = \max_{i=1,\ldots,n-1}\{i : x_i^1 \neq x_i^2\}$. Note that if $x_i^1=x_i^2$ for each $i=1,\ldots, n-1$, then $d({\bf w}_n^1, {\bf w}_n^2)$ can take the values $0,1,2$, depending on the first letters of $w_n^1$ and $w_n^2$ (the vectors correspond to vertices belonging to the same square of side $1$).

By definition of the index $h$, the distance $d(\textbf{w}_n^1,\textbf{ w}_n^2)$ equals the distance $d(\textbf{w}_{h+1}^1,\textbf{w}_{h+1}^2)$ in the graph $\Gamma_{h+1}$; observe that the vectors $\textbf{w}_{h+1}^1$ and $\textbf{w}_{h+1}^2$ are associated with the truncated words $w_{h+1}^1=y^1x_1^1\ldots x_h^1$ and $w_{h+1}^2=y^2x_1^2\ldots x_h^2$, and they occupy two distinct copies of the graph $\Gamma_h$ in $\Gamma_{h+1}$, indexed by $x^1_h$ and $x^2_h$, respectively.\\

CASE I : $d({\bf w}_{h+1}^1, {\bf w}_{h+1}^2)= \|{\bf w}_{h+1}^1 - {\bf w}_{h+1}^2\|_1$, since there is no obstruction in the shortest path from ${\bf w}_{h+1}^1$ to ${\bf w}_{h+1}^2$ in $\Gamma_{h+1}$. This occurs in the following cases:
\begin{itemize}
\item $x^1_h=0, x^2_h =4$; $x^1_h=2, x^2_h =6$.
\item $x^1_h=0,  x^2_h =3$; $x^1_h=1, x^2_h =4$; $x^1_h=2, x^2_h =5$; $x^1_h=3, x^2_h =6$; $x^1_h=4, x^2_h =7$; $x^1_h=5, x^2_h =0$; $x^1_h=6, x^2_h =1$; $x^1_h=7, x^2_h =2$.
\item $x^1_h=1, x^2_h =3$; $x^1_h=3, x^2_h =5$; $x^1_h=5, x^2_h =7$; $x^1_h=7, x^2_h =1$.
\end{itemize}

\indent CASE II\\
$x^1_h=1, x^2_h =5$; we have to consider two different subcases. First of all, observe that the corner vertices of middle hole $H_{h+1}$ in $\Gamma_{h+1}$ are:
$$
A_{h+1}(3^{h-1}, 3^{h-1})       \qquad B_{h+1}(2\cdot 3^{h-1}, 3^{h-1})  \qquad     C_{h+1}(2\cdot 3^{h-1}, 2\cdot 3^{h-1})     \qquad D_{h+1}(3^{h-1}, 2\cdot 3^{h-1}).
$$
Now if $\frac{X({\bf w}_{h+1}^1)+X({\bf w}_{h+1}^2)}{2} \geq \frac{X(A_{h+1})+X(B_{h+1})}{2} = \frac{3^h}{2}$, then:
$$
d({\bf w}_{h+1}^1, {\bf w}_{h+1}^2)= d({\bf w}_{h+1}^1, B_{h+1}) + d(B_{h+1}, C_{h+1}) + d(C_{h+1}, {\bf w}_{h+1}^2) = \|{\bf w}_{h+1}^1 - B_{h+1}\|_1 + 3^{h-1} + \|C_{h+1} - {\bf w}_{h+1}^2\|_1.
$$
Similarly, if $\frac{X({\bf w}_{h+1}^1)+X({\bf w}_{h+1}^2)}{2} < \frac{X(A_{h+1})+X(B_{h+1})}{2}= \frac{3^h}{2}$, then $d({\bf w}_{h+1}^1, {\bf w}_{h+1}^2)= \|{\bf w}_{h+1}^1 - A_{h+1}\|_1 + 3^{h-1} + \|D_{h+1} - {\bf w}_{h+1}^2\|_1.$
An analogous argument holds in the case $x^1_h=3, x^2_h =7$.\\

CASE III\\
$x^1_h=0, x^2_h =2$. Let us put $\ell = \max_{j=1,\ldots, h-1}\{j : x^1_j = 7\}$. We put $\ell = -\infty$ if $x_j^1\neq 7$, for each $j=1,\ldots, h-1$.
\begin{itemize}
\item Case $\ell=-\infty$. We have $d({\bf w}_{h+1}^1, {\bf w}_{h+1}^2) = \|{\bf w}_{h+1}^1 - {\bf w}_{h+1}^2\|_1$, since there is no obstruction in the shortest path from ${\bf w}_{h+1}^1$ to ${\bf w}_{h+1}^2$ in $\Gamma_{h+1}$.
\item Case $\ell\neq -\infty$. Any geodesic path joining ${\bf w}_{h+1}^1$ to ${\bf w}_{h+1}^2$ in $\Gamma_{h+1}$ meets an obstruction (the largest one) given by a hole isomorphic to $H_{\ell +1}$, whose corner vertices are:
\begin{eqnarray}\label{vertices}
A_{\ell+1}=(3^{\ell-1}, 3^{\ell-1}) + \sum_{k=\ell+1}^{h}3^{k-1}v_{x^1_k}  & \qquad & B_{\ell+1}=(2\cdot 3^{\ell-1}, 3^{\ell-1}) + \sum_{k=\ell+1}^{h}3^{k-1}v_{x^1_k} \\ \nonumber
C_{\ell+1}=(2\cdot 3^{\ell-1}, 2\cdot 3^{\ell-1}) + \sum_{k=\ell+1}^{h}3^{k-1}v_{x^1_k} & \qquad & D_{\ell+1}=(3^{\ell-1}, 2\cdot 3^{\ell-1}) + \sum_{k=\ell+1}^{h}3^{k-1}v_{x^1_k}.
\end{eqnarray}
Now if $\frac{Y({\bf w}_{h+1}^1)+Y({\bf w}_{h+1}^2)}{2} \geq \frac{Y(A_{\ell+1})+Y(D_{\ell+1})}{2}$, then:
$$
d({\bf w}_{h+1}^1, {\bf w}_{h+1}^2)= d({\bf w}_{h+1}^1, D_{\ell+1}) + d(D_{\ell+1}, C_{\ell+1}) + d(C_{\ell+1}, {\bf w}_{h+1}^2) = \|{\bf w}_{h+1}^1 - D_{\ell+1}\|_1 + 3^{\ell-1} + \|C_{\ell+1} - {\bf w}_{h+1}^2\|_1.
$$
Similarly, if $\frac{Y({\bf w}_{h+1}^1)+Y({\bf w}_{h+1}^2)}{2} < \frac{Y(A_{\ell+1})+Y(D_{\ell+1})}{2}$, then $d({\bf w}_{h+1}^1, {\bf w}_{h+1}^2)= \|{\bf w}_{h+1}^1 - A_{\ell+1}\|_1 + 3^{\ell-1} + \|B_{\ell+1} - {\bf w}_{h+1}^2\|_1.$

The same argument holds in the case $x^1_h=6, \ x^2_h =4$. Moreover, an analogous method works in the cases $x^1_h=2, x^2_h =4$ and $x^1_h=0, x^2_h =6$, but now the definition of $\ell$ must be replaced with $\ell'=\max_{j=1,\ldots, h-1}\{j : x^1_j = 1\}$, since we have now to consider the obstruction that we meet when we move from the bottom to the top of $\Gamma_{h+1}$.
\end{itemize}

CASE IV\\
$x^1_h=0, x^2_h =1$. Let us define $\ell$ as in Case III.
\begin{itemize}
\item Case $\ell=-\infty$. In this case, we have $d({\bf w}_{h+1}^1, {\bf w}_{h+1}^2)=\|{\bf w}_{h+1}^1 - {\bf w}_{h+1}^2\|_1$, since there is no obstruction in the shortest path from ${\bf w}_{h+1}^1$ to ${\bf w}_{h+1}^2$ in $\Gamma_{h+1}$.\\
\item Case $\ell\neq -\infty$. Any geodesic path connecting ${\bf w}_{h+1}^1$ to ${\bf w}_{h+1}^2$ in $\Gamma_{h+1}$ meets an obstruction (the largest one) given by a hole isomorphic to $H_{\ell +1}$, whose corner vertices are defined as in \eqref{vertices}. If $\frac{Y({\bf w}_{h+1}^1)+Y({\bf w}_{h+1}^2)}{2} \geq \frac{Y(A_{\ell+1})+Y(D_{\ell+1})}{2}$, then:
$$
d({\bf w}_{h+1}^1, {\bf w}_{h+1}^2)= d({\bf w}_{h+1}^1, D_{\ell+1}) + d(D_{\ell+1}, C_{\ell+1}) + d(C_{\ell+1}, {\bf w}_{h+1}^2) = \|{\bf w}_{h+1}^1 - D_{\ell+1}\|_1 + 3^{\ell-1} + \|C_{\ell+1} - {\bf w}_{h+1}^2\|_1.
$$
Similarly, if $\frac{Y({\bf w}_{h+1}^1)+Y({\bf w}_{h+1}^2)}{2} < \frac{Y(A_{\ell+1})+Y(D_{\ell+1})}{2}$, then $d({\bf w}_{h+1}^1, {\bf w}_{h+1}^2)= \|{\bf w}_{h+1}^1 - A_{\ell+1}\|_1 + 3^{\ell-1} + \|B_{\ell+1} - {\bf w}_{h+1}^2\|_1$.
The same argument holds in the cases $x^1_h=1, x^2_h =2$; $x^1_h=6, x^2_h =5$; $x^1_h=5, x^2_h =4$. Finally, a similar argument works in the cases $x^1_h=0, x^2_h =7$; $x^1_h=7, x^2_h =6$; $x^1_h=2, x^2_h =3$; $x^1_h=3, x^2_h =4$, where $\ell$ must be replaced with the index $\ell'$ defined as in the last part of Case III.
\end{itemize}
\begin{example}  \rm
In Fig. \ref{figura} we have represented in $\Gamma_4$ the vertices corresponding to the words $w_4^1=a670$ and $w_4^2=b432$. The corresponding vectors, after the embedding into $\mathbb{Z}^2$, are $(0,5)$ and $(27,5)$, respectively. We have $n=4$, $h=3$, $\ell=2$, so that $A_3(3,3), B_3(6,3), C_3(6,6), D_3(3,6)$ are the corner vertices of the first of the three biggest obstructions met by the shortest path from ${\bf w}_4^1$ to ${\bf w}_4^2$ (case III). One has $d({\bf w}_4^1,{\bf w}_4^2)=29$.

\begin{figure}[h] \psfrag{3}{$3$}\psfrag{4}{$4$}\psfrag{5}{$5$}\psfrag{6}{$6$}\psfrag{7}{$7$} \psfrag{0}{$0$}\psfrag{1}{$1$}\psfrag{2}{$2$}
\psfrag{O}{$O$}\psfrag{w1}{${\bf w}^1_4$}\psfrag{w2}{${\bf w}^2_4$}
\includegraphics[width=0.4\textwidth]{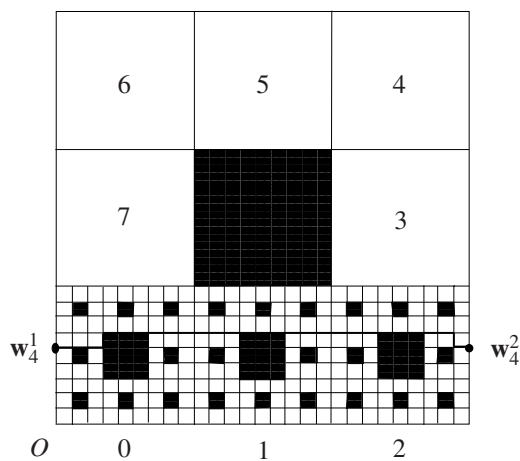}
\caption{Example for the Case III.}\label{figura}
\end{figure}
\end{example}
In the previous description we have redundance if we consider all possible words in $Y\times X^{n-1}$. In fact, as we pointed out before, different words may correspond to the same vertex of $\Gamma_n$. We solve this problem by giving a lexicographic order to such words and considering the smallest one.

\begin{theorem}
The sum of the distances between all vertices obtained in the cases I, II, III, IV and considered without redundance is the Wiener index $W(\Gamma_n)$.
\end{theorem}
The numerical values of $W(\Gamma_n)$ have been computed by using the commercial software Wolfram Mathematica and are reported in Table \ref{tabella}.

\begin{table}[h]
\footnotesize
\begin{tabular}{cccccc}
\hline
$n=1$ & $n=2$ & $n=3$ &  $n=4$ & $n=5$ & $n=6$ \\
\hline 
8 & $320$ & 31264  & 4642456 &    &   \\
\hline 
\end{tabular}
\caption{Wiener index of the carpet graph $\Gamma_n$.} \label{tabella}
\end{table}

\bibliographystyle{aipproc}   

\begin{thebibliography}{99}

\bibitem{bonnier} B.~Bonnier, Y.~Leroyer, and C.~Meyers, Critical exponents for Ising-like systems on
Sierpinski carpets, \emph{J. Physique} {\bf 48},
553--558 (Avril 1987).



\bibitem{rodi} D.~D'Angeli and A.~Donno, Isomorphism classification of infinite Sierpi\'{n}ski carpet graphs, {\it AIP Conference Proceedings} {\bf 1648}, 570002 (2015); doi: 10.1063/1.4912788

\bibitem{compactification} D.~D'Angeli, and A.~Donno, Metric compactification of infinite Sierpi\'{n}ski carpet graphs, preprint, arXiv: 1501.03178


\bibitem{dimeri} D.~D'Angeli, A.~Donno, and T.~Nagnibeda, Counting dimer coverings on self-similar Schreier
graphs, \emph{European J. Combin.} {\bf 33}, Issue 7,
1484--1513 (2012).

\bibitem{ising} D.~D'Angeli, A.~Donno, and T.~Nagnibeda, Partition functions of the Ising model on some self-similar Schreier
graphs, in \emph{Progress in Probability: Random Walks, Boundaries
and Spectra} {\bf 64}, edited by D.~Lenz, F.~Sobieczky and W.~Woess, Springer Basel (2011), pp. 277--304.

\bibitem{sandi1} A.~A.~Dobrynin, R.~Entringer, and I.~Gutman, Wiener index of trees: Theory and applications,
\textit{Acta Appl. Math.} \textbf{66}, 211--249 (2001).

\bibitem{sandi2} A.~A.~Dobrynin, I.~Gutman, S.~Klav\v{z}ar, and P.~\v{Z}igert, Wiener index of hexagonal systems,
\textit{Acta Appl. Math.} \textbf{72}, 247--294 (2002).

\bibitem{donnoiacono} A.~Donno and D.~Iacono, Distances and isomorphisms in $4$-regular circulant graphs, \emph{Proceedings of the 2nd Minisymposium on Mathematics in Engineering and Technology}, ICNAAM 2015, Rhodes, 23--29/09/2015, accepted.


\bibitem{gefen2} Y.~Gefen, A.~Aharony, Y.~Shapir, and B.~B.~Mandelbrot, Phase transitions on fractals. II. Sierpi\'{n}ski gaskets, \emph{J. Phys. A} {\bf 17}, no. 2, 435--444  (1984).

\bibitem{gefen3} Y.~Gefen, A.~Aharony, and B.~B.~Mandelbrot, Phase
transitions on fractals. III. Infinitely ramified lattices,
\emph{J. Phys. A} {\bf 17}, no. 6, 1277--1289  (1984).


\bibitem{percolation} M.~Shinoda, Existence of phase transition of percolation on Sierpi\'{n}ski
carpet lattices, \emph{J. Appl. Probab.} {\bf 39}, no. 1,
1--10  (2002).

\bibitem{courbe} W.~Sierpi\'{n}ski, Sur une courbe cantorienne qui contient une image
biunivoque et continue de toute courbe donn\'{e}e, \emph{C. R. Acad.
Sci. Paris}, \textbf{162}, 629--642  (1916).


\bibitem{Wiener} H.~Wiener, Structural determination of paraffin boiling points, \textit{J. Amer. Chem. Soc.} \textbf{69} 17--20 (1947).

\end{thebibliography}

 \end{document}